\documentclass[12pt,a4paper,reqno]{amsart}
\usepackage{amsmath}
\usepackage{amsfonts}
\usepackage{amssymb}
\numberwithin{equation}{section}

\theoremstyle{plain}

\newtheorem{theorem}[subsection]{Theorem}
\newtheorem{proposition}[subsection]{Proposition}
\newtheorem{lemma}[subsection]{Lemma}

\theoremstyle{definition}

\newtheorem{definition}[subsection]{Definition}

\renewcommand{\leq}{\leqslant}
\renewcommand{\geq}{\geqslant}

\newsavebox{\proofbox}
\savebox{\proofbox}{\begin{picture}(7,7)%
  \put(0,0){\framebox(7,7){}}\end{picture}}

\newcommand{\md}[1]{\ensuremath{(\mbox{mod}\, #1)}}

\def\Z{{\mathbb Z}}
\def\E{{\mathbb E}}

\def\R{{\mathbb R}}

\def\Q{{\mathbb Q}}

\def\x{{\bf x}}

\def\w{w}

\def\vs{\vspace{11pt}}
\def\ni{\noindent}

\def\proof{\noindent\textit{Proof. }}
\def\endproof{\hfill{\usebox{\proofbox}}}

\def\emph#1{{\it #1}}
\def\textbf#1{{\bf #1}}

\begin{document}

\title{Arithmetic progressions of primes in short intervals}

\author{Chunlei Liu}
\address{School of Mathematical Science\\
Beijing Normal University\\
Beijing 100875\\
} \email{clliu@bnu.edu.cn}

\thanks{This work is supported by Project 10671015 of the Natural Science Foundation of China.}


\begin{abstract}
\ni Green and Tao proved that the primes contains arbitrarily long
arithmetic progressions. We show that, essentially the same proof
leads to the following result: If $N$ is sufficiently large and $M$
is not {\it too} small compared with $N$, then the primes in the
interval $[N,N+M]$ contains {\it many} arithmetic progressions of
length $k$.
\end{abstract}

\maketitle

\section{Introduction}

\ni Let $N$ be a positive integer going to infinity. We write $o(1)$
for any quantity which tends to zero as $N$ goes to infinity, and
write $O(1)$ for any quantity which has a bound independent of $N$.
Let $w=w(N)\leq\frac12\log\log N$ be any function which tends to
infinity with $N$, and let $W: = \prod_{p \leq w} p$ be the product
of the primes up to $w$. Let $\tilde \Lambda$ be the $W$-tricked von
Mangoldt function defined by
\[ \widetilde{\Lambda}(n) := \left\{
\begin{array}{ll}
\frac{\phi(W)}{W} \log(Wn + 1) & \hbox{ when } Wn+1 \hbox{ is prime}\\
0 & \hbox{ otherwise.}
\end{array}\right.\]
\vs \ni Let $M$ be a large prime number. Define $\Z_M:=\Z/M\Z$ to be
the finite field consisting of residue classes modulo $M$. We always
identify $\Z_M$ with the set
$$\{N-M,N-M+1,\cdots,N-1\},$$ which is a complete system of
representatives modulo $M$.\vs

\ni If $A$ is a finite non-empty set and $f: A \to \R$ is a
function, we write $$\E( f ) := \E( f(x) | x \in A )$$ for the
average value of $f$, that is to say
\[ \E(f) := \frac{1}{|A|} \sum_{x \in A} f(x).\] Here, as is usual, we write $|A|$ for the cardinality of the set $A$.
More generally, if $P(x)$ is any statement concerning an element of $A$ which is true for at least one $x \in A$,
we define
\[ \E(f(x) | P(x)) := \frac{ \sum_{x \in A: P(x)} f(x)}{|\{ x \in A: P(x) \}|}.\]
\vs

\ni Let $k$ be any fixed integer greater than $3$, and let
$\epsilon_k := 1/2^{k}(k+4)!$. \vs

\ni A famous theorem of Green-Tao in \cite{green-tao} asserts that
the prime numbers contain arbitrarily long arithmetic progressions.
In this paper we show that the proof of Green-Tao really yields the
following theorem.
\begin{theorem}\label{primes-short-interval}Let $M$ be a function of $N$ with values in the set of prime
numbers which satisfies $N^{\varepsilon}<M\leq N$ for some positive
number $\varepsilon$. Suppose that on the interval
$[N+\epsilon_kM,N+2\epsilon_kM]$ the mean value of the $W$-tricked
von Mongoldt function tends to $1$ as $N$ goes to infinity. Define
the function $f$ on $\Z_M$ by setting
$$f(n) :=\left\{
\begin{array}{ll}
k^{-1}2^{-k-5}\tilde{\Lambda}(n)&\hbox{ when }\epsilon_k M \leq n-N
\leq 2\epsilon_k M\\
0 &\hbox{ otherwise.}\end{array}\right.$$ Then there is a positive
constant $c_k$ depending only on $k$ such that
$$\mathbb{E} \big( f(x)f(x+r) \dots f(x + (k-1)r) \; \big| \; x,r \in \mathbb{Z}_M\big) \geq c_k- o(1).$$
\end{theorem}\vs

\ni From that theorem we see that, for sufficiently large $N$, there
are at least $b_kM^2/\log^kN$ arithmetic progressions of length $k$
consisting of primes in the interval $(WN,W(N+M)]$, where $b_k$ is a
positive constant $c_k$ depending only on $k$. According to
Green-Tao, we can in fact take $w$ to be a sufficiently large number
independent of $N$, depending only on $k$. Then $W$ will be a
constant depending only on $k$.\vs

\ni\textit{Acknowledgements}The author would like to thank the
Morningside Center of Mathematics, Chinese Academy of Sciences for
support over several years.

\section{The linear forms property}

\ni In this section we construct a majorant $\nu$ for $f$ and prove
that $\nu$ satisfies the linear forms condition.

\begin{definition}\label{truncated-divisor-sum} Let $R$ be a parameter (in applications it will be a small power of $N$). Define
\[ \Lambda_R(n) := \sum_{\substack{d | n \\ d \leq R}} \mu(d)\log (R/d)
= \sum_{d|n} \mu(d) \log(R/d)_+.\]
\end{definition}

\ni These truncated divisor sums have been studied in several
papers, most notably the works of Goldston and Y{\i}ld{\i}r{\i}m
\cite{goldston-yildirim-old1,goldston-yildirim-old2,goldston-yildirim}
concerning the problem of finding small gaps between primes.

\begin{definition}\label{mu-def} Let $R := M^{k^{-1}2^{-k-4}}$.  We define the function
$\nu: \Z_M \to \R^+$ by
\[ \nu(n) \; := \; \left\{
\begin{array}{ll}
\frac{\phi(W)}{W} \frac{\Lambda_R(Wn + 1)^2}{\log R} & \hbox{ when } \epsilon_k M \leq n-N \leq 2\epsilon_k M\\
1 & \hbox{ otherwise}
\end{array}\right.
\]
for all $N+M \leq n < N+M$.
\end{definition}

\begin{lemma}\label{mu-majorises}
Let $N$ be a sufficiently large integer depending on $k$. Then the
function $\nu$ is majorant for $f$ in Theorem
\ref{primes-short-interval}. That is, $\nu(n) \geq 0$ for all $n \in
\mathbb{Z}_M$, and $\nu(n) \geq
k^{-1}2^{-k-5}\widetilde{\Lambda}(n)$ for all $N+\epsilon_k M \leq n
\leq N+2\epsilon_kM$.
\end{lemma}

\proof The first claim is trivial. The second claim is also trivial
unless $Wn+1$ is prime.  From definition of $R$, we see that $Wn+1 >
R$ if $N$ is sufficiently large.  Then the sum over $d | Wn + 1$, $d
\leq R$ in \eqref{truncated-divisor-sum} in fact consists of just
the one term $d = 1$. Therefore $\Lambda_R(Wn + 1) = \log R$, which
means that $\nu(n) = \frac{\phi(W)}{W}\log R \geq k^{-1}2^{-k-5}
\widetilde{\Lambda}(n)$ by construction of $R$ and $N$.
\endproof

\vspace{11pt}\begin{definition}[Linear forms
condition]\label{linear-forms-condition} Let $m_0,t_0$ and $L_0$ be
small positive integer parameters. Then we say that $\nu:
\mathbb{Z}_M \rightarrow \mathbb{R}^{+}$ satisfies the
$(m_0,t_0,L_0)$-linear forms condition if the following holds. Let
$m \leq m_0$ and $t \leq t_0$ be arbitrary, and suppose that $(
L_{ij} )_{1 \leq i \leq m, 1 \leq j \leq t}$ are arbitrary rational
numbers with numerator and denominator at most $L_0$ in absolute
value, and that $b_i$, $1 \leq i \leq m$, are arbitrary elements of
$\mathbb{Z}_M$.  For $1 \leq i \leq m$, let $\psi_i: \Z_M^t \to
\Z_M$ be the linear forms $\psi_i(\mathbf{x}) = \sum_{j = 1}^t
L_{ij}x_j + b_i$, where $\mathbf{x} = (x_1,\dots,x_t) \in
\mathbb{Z}_M^t$, and where the rational numbers $L_{ij}$ are
interpreted as elements of $\Z_M$ in the usual manner (assuming $M$
is prime and larger than $L_0$). Suppose that as $i$ ranges over
$1,\ldots,m$, the $t$-tuples $(L_{ij})_{1 \leq j \leq t} \in \Q^t$
are non-zero, and no $t$-tuple is a rational multiple of any other.
Then we have
\begin{equation}\label{lfc}
 \mathbb{E}\left( \nu(\psi_1(\mathbf{x})) \dots \nu(\psi_m(\mathbf{x})) \;|\; \mathbf{x} \in \mathbb{Z}_M^t \right) = 1 + o_{L_0,m_0,t_0}(1).
\end{equation}
Note that the rate of decay in the $o(1)$ term is assumed to be
uniform in the choice of $b_1, \ldots, b_m$.
\end{definition}

\ni The following propositions plays a crucial role in proving that
$\nu$ satisfies the linear forms condition.

\begin{proposition}[Goldston-Y{\i}ld{\i}r{\i}m]\label{GY}  Let $m, t$ be positive integers.  For each $1 \leq i \leq m$, let
$\psi_i(\mathbf{x}) := \sum_{j=1}^t L_{ij} x_j + b_i$, be linear
forms with integer coefficients $L_{ij}$ such that $|L_{ij}| \leq
\sqrt{\w(N)}/2$ for all $i = 1,\ldots m$ and $j=1,\ldots,t$. We
assume that the $t$-tuples $(L_{ij})_{j=1}^t$ are never identically
zero, and that no two $t$-tuples are rational multiples of each
other. Write $\theta_i := W\psi_i + 1$.  Suppose that $B$ is a
product $\prod_{i = 1}^t I_i \subset \R^t$ of $t$ intervals $I_i$,
each of which having length  at least $R^{10m}$. Then
\[ \E( \Lambda_R(\theta_1(\mathbf{x}))^2 \dots \Lambda_R(\theta_m(\mathbf{x}))^2 | \x \in B ) =
(1 + o_{m,t}(1)) \left(\frac{W \log R}{\phi(W)}\right)^m.
\]
\end{proposition}
\noindent\textit{Remarks.} That proposition was stated and proved by
Green-Tao in \cite{green-tao}, however, according to Green-Tao, it
is a straightforward generalisation of \cite[Proposition
2]{goldston-yildirim}.\vs

\ni Before proving the linear forms condition, We show that
$\mathbb{E}(\nu) = 1 + o(1)$.

\begin{lemma}\label{nu-measure} We have
$\mathbb{E}(\nu) = 1 + o(1)$.
\end{lemma}
\proof Apply Proposition \ref{GY} with $m := t := 1$, $\psi_1(x_1)
:= x_1$ and $B := [N+\epsilon_k M, N+2\epsilon_k M]$ (taking $N$
sufficiently large depending on $k$, of course). Comparing with
Definition \ref{mu-def} we thus have
\[ \E (\nu(x) \; | \; x \in [N+\epsilon_k M, N+2\epsilon_k M]) = 1 + o(1).\]
But from the same definition we clearly have
\[ \E (\nu(x) \; | \; x \in \Z_M \backslash [N+\epsilon_k M, N+2\epsilon_k M] ) = 1;\]
Combining these two results confirms the lemma.\endproof\vs

\ni Now we verify the linear forms condition, which is proven in a
similar spirit to the above lemma.

\begin{proposition}\label{mu-linear-forms}
The function $\nu$ satisfies the $(k \cdot 2^{k-1},3k-4,k)$-linear
forms condition.
\end{proposition}
\proof Let $\psi_i(x) = \sum_{j = 1}^t L_{ij} x_j + b_i$ be linear
forms of the type which feature in Definition
\ref{linear-forms-condition}. That is to say, we have $m \leq k
\cdot 2^{k-1}$, $t \leq 3k - 4$, the $L_{ij}$ are rational numbers
with numerator and denominator at most $k$ in absolute value, and
none of the $t$-tuples $(L_{ij})_{j=1}^t$ is zero or is equal to a
rational multiple of any other. We wish to show that
\begin{equation}\label{to-prove} \mathbb{E}(\nu(\psi_1(\mathbf{x}))
\dots \nu(\psi_m(\mathbf{x})) \; | \; \mathbf{x} \in \mathbb{Z}_M^m)
= 1 + o(1).\end{equation}

\ni We may clear denominators and assume that all the $L_{ij}$ are
integers, at the expense of increasing the bound on $L_{ij}$ to
$|L_{ij}| \leq (k+1)!$.  Since $\w(N)$ is growing to infinity in
$N$, we may assume that $(k+1)! < \sqrt{\w(N)}/2$ by taking $N$
sufficiently large. This is required in order to apply Proposition
\ref{GY} as we have stated it. \vs

\ni The two-piece definition of $\nu$ in Definition \ref{mu-def}
means that we cannot apply Proposition \ref{GY} immediately, and we
need the following localization argument. \vs

\ni We chop the range of summation in \eqref{to-prove} into $Q^t$
almost equal-sized boxes, where $Q = Q(N)$ is a slowly growing
function of $N$ to be chosen later. Thus let
\[ B_{u_1,\dots,u_t} = \{ \mathbf{x} \in \Z_M^m : x_j \in [N+\lfloor u_jM/Q\rfloor,N+\lfloor(u_j + 1)M/Q\rfloor), j = 1,\dots,t\},\] where the $u_j$ are to be considered $\md{Q}$. Observe that up to negligible multiplicative errors of $1 + o(1)$ (arising because the boxes do not quite have equal sizes) the left-hand side of \eqref{to-prove} can be rewritten as
$$ \E( \E( \nu(\psi_1(\x)) \ldots \nu(\psi_m(\x)) | \x \in B_{u_1,\ldots,u_t} ) | u_1, \ldots, u_t \in \Z_Q ).$$
Call a $t$-tuple $(u_1,\dots,u_t) \in \Z_Q^t$ \textit{nice} if for
every $1 \leq i \leq m$, the sets $\psi_i(B_{u_1,\ldots,u_t})$ are
either completely contained in the interval $[N+\epsilon_k M,
N+2\epsilon_k M]$ or are completely disjoint from this interval.
From Proposition \ref{GY} and Definition \ref{mu-def} we observe
that
$$ \E( \nu(\psi_1(\x)) \ldots \nu(\psi_m(\x)) | \x \in B_{u_1,\ldots,u_t} ) = 1 + o_{m,t}(1)$$
whenever $(u_1,\ldots,u_t)$ is nice, since we can replace
 each of the $\nu(\psi_i(\x))$ factors by either
$\frac{\phi(W)}{W \log R} \Lambda_R^2(\theta_i(\x))$ or $1$, and
$M/Q$ will exceed $R^{10m}$ for $Q$ sufficiently slowly growing in
$N$, by definition of $R$ and the upper bound on $m$. When
$(u_1,\ldots,u_t)$ is not nice, then we can crudely bound $\nu$ by
$1 + \frac{\phi(W)}{W \log R} \Lambda_R^2(\theta_i(\x))$, multiply
out, and apply Proposition \ref{GY} again to obtain
$$ \E( \nu(\psi_1(\x)) \ldots \nu(\psi_m(\x)) | \x \in B_{u_1,\ldots,u_t} ) = O_{m,t}(1) + o_{m,t}(1)$$
We shall shortly show that the proportion of non-nice $t$-tuples
$(u_1,\ldots,u_t)$ in $\Z_Q^t$ is at most $O_{m,t}(1/Q)$, and thus
the left-hand side of \eqref{to-prove} is $1 + o_{m,t}(1) +
O_{m,t}(1/Q)$, and the claim follows by choosing $Q$ sufficiently
slowly growing in $N$.\vs

\ni It remains to verify the claim about the proportion of non-nice
$t$-tuples. Suppose $(u_1,\ldots,u_t)$ is not nice. Then there
exists $1 \leq i \leq m$ and $\x, \x' \in B_{u_1,\ldots,u_t}$ such
that $\psi_i(\x)$ lies in the interval $[N+\epsilon_k M,
N+2\epsilon_k M]$, but $\psi_i(\x')$ does not.  But from definition
of $B_{u_1,\ldots,u_t}$(and the boundedness of the $L_{ij}$)  we
have
$$ \psi_i(\x), \psi_i(\x') = \sum_{j=1}^t L_{ij}(N+\lfloor Mu_j/Q\rfloor) + b_i + O_{m,t}(M/Q). $$
Thus we must have
$$ N+a \epsilon_k M = \sum_{j=1}^t L_{ij}(N+\lfloor Mu_j/Q\rfloor) + b_i + O_{m,t}(M/Q)$$
for either $a=1$ or $a=2$.  Dividing by $M/Q$, we obtain
$$ \sum_{j=1}^t L_{ij} u_j =(1-\sum_{j=1}^t L_{ij}-b_i)Q/N +a \epsilon_k Q + O_{m,t}(1) \quad \md{Q}.$$
Since $(L_{ij})_{j=1}^t$ is non-zero, the number of $t$-tuples
$(u_1,\ldots,u_t)$ which satisfy this equation is at most
$O_{m,t}(Q^{t-1})$.  Letting $a$ and $i$ vary we thus see that the
proportion of non-nice $t$-tuples is at most $O_{m,t}(1/Q)$ as
desired (the $m$ and $t$ dependence is irrelevant since both are
functions of $k$).
\endproof\vspace{11pt}
\section{The correlation property}
\ni In this section we show that $\nu$ satisfies the correlation
condition.
\begin{definition}[Correlation condition]\label{correlation-condition}
Let $m_0$ be a positive integer parameter. We say that $\nu:
\mathbb{Z}_M \rightarrow \mathbb{R}^{+}$ satisfies the
$m_0$-correlation condition if for every $1 < m \leq m_0$ there
exists a weight function $\tau = \tau_{m}: \mathbb{Z}_M \to \R^+$
which obeys the moment conditions
\begin{equation}\label{eq3.1}
 \E( \tau^q ) = O_{m,q}(1)
\end{equation}
for all $1 \leq q < \infty$ and such that
\begin{equation}\label{eq3.2}
\E( \nu(x+h_1) \nu(x+h_2) \ldots \nu(x+h_m) \; |\;  x \in
\mathbb{Z}_M) \leq \sum_{1 \leq i < j \leq m} \tau(h_i-h_j)
\end{equation}
for all $h_1, \ldots, h_m \in \Z_M$ (not necessarily distinct).
\end{definition}

\ni The following proposition plays a crucial roles in proving that
$\nu$ satisfies the correlation condition.
\begin{proposition}[Goldston-Y{\i}ld{\i}r{\i}m]\label{GY2}
Let $m \geq 1$ be an integer, and let $B$ be an interval of length
at least $R^{10m}$. Suppose that $h_1,\dots,h_m$ are distinct
integers satisfying $|h_i| \leq N^2$ for all $1 \leq i \leq m$, and
let $\Delta$ denote the integer
\[ \Delta := \prod_{1 \leq i < j \leq m} |h_i - h_j|.\]
Then
\begin{equation}\label{GY2-est}
\begin{split}
 \E( &\Lambda_R(W(x_1 + h_1) + 1)^2 \dots \Lambda_R(W(x_m + h_m) + 1)^2| x \in B ) \\
&\leq (1 + o_m(1)) \left( \frac{W\log R}{\phi(W)}\right)^m \prod_{p
| \Delta}(1 + O_m(p^{-1/2})).
\end{split}
\end{equation}
Here and in the sequel, $p$ is always understood to be prime.
\end{proposition}

\noindent\textit{Remarks.} That proposition was stated and proved by
Green-Tao in \cite{green-tao}, however, Green-Tao attributed it to
Goldston-Y{\i}ld{\i}r{\i}m for reasons similar to Proposition
\ref{GY}.\vs

\ni In a short while we will use Proposition \ref{GY2} to show that
$\nu$ satisfies the correlation condition. Prior to that, however,
we must look at the average size of the ``arithmetic'' factor
$\prod_{p | \Delta} (1 + O_m(p^{-1/2}))$ appearing in that
proposition.

\begin{lemma}\label{additive-weights}  Let $m \geq 1$ be a parameter.
There is a weight function $\tau = \tau_m: \Z \to \R^+$ such that
$\tau(n) \geq 1$ for all $n \neq 0$, and such that for all distinct
$h_1,\ldots,h_j \in [N+\epsilon_k M, N+2\epsilon_k M]$ we have
\[ \prod_{p | \Delta}(1 + O_m(p^{-1/2})) \leq \sum_{1 \leq i < j \leq m} \tau(h_i - h_j),\]
where $\Delta$ is defined in Proposition \ref{GY2}, and such that
$\E( \tau^q(n) | 0 < |n| \leq M ) = O_{m,q}(1)$ for all $0< q <
\infty$.
\end{lemma}

\proof We observe that
$$ \prod_{p|\Delta} (1 + O_m(p^{-1/2})) \leq \prod_{1 \leq i < j \leq m}
\bigg(\prod_{p|h_i-h_j} (1 + p^{-1/2}) \bigg)^{O_m(1)}.$$
By the arithmetic mean-geometric mean inequality (absorbing all
constants into the $O_m(1)$ factor) we can thus take $\tau_m(n) :=
O_m(1) \prod_{p | n} (1 + p^{-1/2})^{O_m(1)}$ for all $n \neq 0$.
(The value of $\tau$ at 0 is irrelevant for this lemma since we are
taking all the $h_i$ to be distinct). To prove the claim, it thus
suffices to show that
$$ \E\bigg( \prod_{p|n} (1 + p^{-1/2})^{O_m(q)} \; \bigg| \; 0 < |n| \leq M \bigg)
= O_{m,q}(1) \hbox{ for all } 0 < q < \infty.$$
Since $(1 + p^{-1/2})^{O_m(q)}$ is bounded by $1 + p^{-1/4}$ for all
but $O_{m,q}(1)$ many primes $p$, we have
$$ \E\bigg( \prod_{p|n} (1 + p^{-1/2})^{O_m(q)} \;  \bigg| \; 0 < |n| \leq M \bigg) \leq O_{m,q}(1)
\E\bigg( \prod_{p|n} (1 + p^{-1/4}) \; \bigg| \; 0 < n \leq M
\bigg).$$ But $\prod_{p|n} (1 + p^{-1/4}) \leq \sum_{d|n} d^{-1/4}$,
and hence
\begin{eqnarray*} \E\bigg( \prod_{p|n} (1 + p^{-1/2})^{O_m(q)} \; \bigg| \; 0 < |n| \leq M \bigg) & \leq  & O_{m,q}(1)
\frac{1}{2M} \sum_{1 \leq |n| \leq M} \sum_{d|n} d^{-1/4} \\ & \leq
& O_{m,q}(1) \frac{1}{2M} \sum_{d=1}^M \frac{M}{d}
d^{-1/4},\end{eqnarray*} which is $O_{m,q}(1)$ as
desired.\endproof\vspace{11pt}

\ni We are now ready to verify the correlation condition.

\begin{proposition}\label{mu-cor-con} The measure $\nu$ satisfies the $2^{k-1}$-correlation condition.\end{proposition}
\proof Let us begin by recalling what it is we wish to prove. For
any $1 \leq m \leq 2^{k-1}$ and $h_1,\dots,h_m \in \mathbb{Z}_N$ we
must show a bound
\begin{equation}\label{to-prove2} \E\big( \nu(x+h_1) \nu(x+h_2) \ldots \nu(x+h_m) \; \big| \; x \in \mathbb{Z}_N\big)
\leq \sum_{1 \leq i < j \leq m} \tau(h_i-h_j),\end{equation} where
the weight function $\tau= \tau_m$ is bounded in $L^q$ for all
$q$.\vs

\ni Fix $m$, $h_1, \ldots, h_m$. We shall take the weight function
constructed in Lemma \ref{additive-weights} (identifying $\Z_M$ with
the integers between $-M/2$ and $+M/2$), and set \[ \tau(0) :=
\exp(Cm \log N/\log \log N)\] for some large absolute constant $C$.
From the previous lemma we see that $\E(\tau^q) = O_{m,q}(1)$ for
all $q$, since the addition of the weight $\tau(0)$ at 0 only
contributes $o_{m,q}(1)$ at most.\vs

\ni We first dispose of the easy case when at least two of the $h_i$
are equal.  In this case we bound the left-hand side of
\eqref{to-prove} crudely by $\| \nu \|_{L^\infty}^m$.  But from
Definitions \ref{truncated-divisor-sum}, \ref{mu-def} and by
standard estimates for the maximal order of the divisor function
$d(n)$ we have the crude bound $\| \nu \|_{L^\infty} \ll \exp(C\log
N/\log \log N)$, and the claim follows thanks to our choice of
$\tau(0)$.\vs

\ni Suppose then that the $h_i$ are distinct. Write
$$g(n) := \frac{\phi(W)}{W} \frac{\Lambda_R^2(Wn+1)}{\log R} {\bf 1}_{[N+\epsilon_k M, N+2\epsilon_k M]}(n).$$
Then by construction of $\nu$ (Definition \ref{mu-def}), we have
\begin{eqnarray*} &&\mathbb{E}\big(\nu(x + h_1) \dots \nu(x+h_m) \; \big| \; x \in \mathbb{Z}_M\big) \\
&& \qquad\qquad\qquad\leq  \mathbb{E}\big((1 + g(x + h_1)) \dots (1
+ g(x+h_m)) \; \big| \; x \in \mathbb{Z}_M\big).\end{eqnarray*} The
right-hand side may be rewritten as
\[ \sum_{A \subseteq \{1,\ldots,m\}}  \mathbb{E}\bigg( \prod_{i \in A} g(x+h_i) \; \bigg| \; x \in \mathbb{Z}_M\bigg)\]
 Observe that for $i,j \in A$ we may assume $|h_i - h_j| \leq
\epsilon_k M$, since the expectation vanishes otherwise.  By
Proposition \ref{GY2} and Lemma \ref{additive-weights}, we therefore
have
\[ \mathbb{E}\bigg( \prod_{i \in A} g(x+h_i) \; \bigg| \; x \in \mathbb{Z}_M\bigg) \leq \sum_{1 \leq i < j \leq m} \tau(h_i - h_j) + o_m(1).\]
Summing over all $A$, and adjusting the weights $\tau$ by a bounded
factor (depending only on $m$ and hence on $k$), we obtain the
result.
\endproof\vspace{11pt}

\section{Proof of the main theorem}

\ni In this section we conclude the proof of Theorem
\ref{primes-short-interval}.

\begin{definition}\label{pseudo-def}
Let $\nu : \mathbb{Z}_M \rightarrow \mathbb{R}^{+}$ be a function.
We say that $\nu$ is $k$-pseudorandom measure if it obeys the
estimate $\mathbb{E}(\nu)=1+o(1)$ and satisfies the $(k\cdot
2^{k-1},3k - 4,k)$-linear forms condition as well as the
$2^{k-1}$-correlation condition.
\end{definition}

\begin{theorem}[Green-Tao]\label{prime-majorant} The function $\nu :
\mathbb{Z}_M \rightarrow \R^+$ in Definition \ref{mu-def} is a
$k$-pseudorandom measure that majorises $f$ in Theorem
\ref{primes-short-interval}
\end{theorem}
\proof That theorem follows from Lemmas \ref{mu-majorises},
\ref{nu-measure} and Propositions \ref{mu-linear-forms},
\ref{mu-cor-con}.\endproof\vs

\noindent\textit{Remarks.} I have attributed this theorem to Green
and Tao, because the above argument is a straightforward
generalisation of that of \cite[Proposition 9.1]{green-tao}.\vs

\ni The proof of Theorem \ref{primes-short-interval} is base on the
following theorem.

\begin{theorem}[Green-Tao]\label{green-tao-sz} Let $k \geq 3$ and $0 < \delta
\leq 1$  be fixed parameters. Suppose that $\nu: \mathbb{Z}_M \to
\R^+$ is $k$-pseudorandom measure. Let $f: \mathbb{Z}_M \to \R^+$ be
any non-negative function obeying the bound
\begin{equation}\label{f-bound-2}
0 \leq f(x) \leq \nu(x) \hbox{ for all } x \in \mathbb{Z}_M
\end{equation}
and
\begin{equation}\label{f-density}
 \E(f) \geq \delta.
\end{equation}
Then we have
\begin{equation}\label{recurrence}
 \E( f(x) f(x+r) \ldots f(x+(k-1)r) | x,r \in \mathbb{Z}_M ) \geq c(k,\delta) - o_{k,\delta}(1)
\end{equation}
where $c(k,\delta) > 0$ stands for a constant depending only on $k$
and $\delta$.
\end{theorem}
\vs

\ni That theorem is a great generalization of the following theorem.
\begin{theorem}[Szemer\'edi's theorem]\label{sz}Let $k
\geq 3$ and $0 < \delta \leq 1$  be fixed parameters. Let $f:
\mathbb{Z}_M \to \R^+$ be any function which is bounded by a bound
independent of $M$. Suppose that
\begin{equation}
 \E(f) \geq \delta.
\end{equation}
Then we have
\begin{equation}\label{recurrence}
 \E( f(x) f(x+r) \ldots f(x+(k-1)r) | x,r \in \mathbb{Z}_M ) \geq c(k,\delta) - o_{k,\delta}(1)
\end{equation}
where $c(k,\delta) > 0$ is the same constant which appears in
Theorem \ref{green-tao-sz}. \textup{(}The decay rate
$o_{k,\delta}(1)$, on the other hand, decays significantly faster
than that in Theorem \ref{green-tao-sz}\textup{)}.
\end{theorem}
\vs

\noindent\textit{Remarks.} The $k=3$ case of Szemer\'edi's theorem
was established by Roth\cite{roth}. The general case as well as the
$k=4$ case was proved by Szemer\'edi \cite{szemeredi-4,szemeredi}.
The formulation here is different from the original one, but can be
deduce from the original one. The argument was first worked out by
Varnavides \cite{varnavides}). A direct proof of Theorem \ref{sz}
can be found in \cite{tao:ergodic}.\vs

\noindent\textit{Proof of Theorem \ref{primes-short-interval}.} By
our assumption on $M$, we see that
$$ \E(f) = \frac{k^{-1}2^{-k-5}}{M} \sum_{N+\epsilon_k M \leq n \leq N+2 \epsilon_k M} \tilde \Lambda(n)
= k^{-1}2^{-k-5} \epsilon_k (1 + o(1)).$$ We now apply Theorem
\ref{prime-majorant} and Theorem \ref{green-tao-sz} to conclude that
$$\mathbb{E} \big( f(x)f(x+r) \dots f(x + (k-1)r) \; \big| \; x,r \in \mathbb{Z}_N\big) \geq c(k,k^{-1}2^{-k-5} \epsilon_k) - o(1).$$
Theorem \ref{primes-short-interval} follows by setting
$c_k=c(k,k^{-1}2^{-k-5} \epsilon_k)$.
\endproof\vspace{11pt}

\providecommand{\bysame}{\leavevmode\hbox to3em{\hrulefill}\thinspace}
\providecommand{\MR}{\relax\ifhmode\unskip\space\fi MR }
\providecommand{\MRhref}[2]{%
  \href{http://www.ams.org/mathscinet-getitem?mr=#1}{#2}
}
\providecommand{\href}[2]{#2}

     \end{document}